\magnification 1200
\def\R{{\rm I\kern-0.2em R\kern0.2em \kern-0.2em}}
\def\N{{\rm I\kern-0.2em N\kern0.2em \kern-0.2em}}
\def\P{{\rm I\kern-0.2em P\kern0.2em \kern-0.2em}}
\def\B{{\rm I\kern-0.2em B\kern0.2em \kern-0.2em}}
\def\Z{{\rm I\kern-0.2em Z\kern0.2em \kern-0.2em}}
\def\C{{\bf \rm C}\kern-.4em {\vrule height1.4ex width.08em depth-.04ex}\;}
\def\B{{\bf \rm B}\kern-.4em {\vrule height1.4ex width.08em depth-.04ex}\;}

\def\D{{\Delta}}

\def\bD{{b\Delta}}

\def\cC{{\cal C}}
\def\cP{{\cal P}}

\def\cC{{\cal C}}
\def\cS{{\cal S}}
\def\cS{{\cal S}}

\def\cP{{\cal P}}

\font\ninerm=cmr8
\
\vskip 25mm
\centerline {\bf THE WINDING NUMBER OF PF+1 FOR POLYNOMIALS P}

\centerline {\bf AND MEROMORPHIC EXTENDIBILITY OF F}
\vskip 4mm
\centerline{Josip Globevnik}
\vskip 4mm
{\noindent \ninerm ABSTRACT\ \ Let $\D$ be the open unit disc in $\C$. 
The paper deals with the following conjecture: \ If $f$ is a continuous function on 
$\bD$ such that the change of argument of $Pf+1$ around $b\D$  is nonnegative 
for every polynomial $P$ such that $Pf+1$ has no zero on $b\D$ 
then $f$ extends holomorphically through $\D$. We prove a related result on 
meromorphic extendibility for smooth functions with finitely many zeros of 
finite order, which, in particular, implies that   
the conjecture holds for real analytic functions. 
} 
\vskip 4mm
\bf 1.\ Introduction \rm

Let $\D$ be the open unit disc in $\C$. Given a continuous function $\varphi $ on $b\D$ with no 
zero on $b\D$ we denote by $W(\varphi )$ the winding number of $\varphi $ around $0$ so that 
$2\pi W(\varphi)$ is the change of argument of $\varphi (z)$ as $z$ runs around $b\D $ in positive direction.  
If a function $g$ is holomorphic on $\D $ then we 
denote by $Z(g)$ the number of zeros of $g$ counting multiplicity. We denote by $A( \D )$ the disc algebra,
 that is the algebra of all continuous functions on $\overline \D$ which are holomorphic on $\D$. 
It is known that one can characterize holomorphic extendibility in terms of the argument principle:
\vskip 2mm
\noindent \bf THEOREM 1.1\ [G1, G2, K] \ \it A continuous function $f$ on $b\D$ extends holomorphically through $\D$ if and only if $W(f+Q)\geq 0$ 
for every polynomial $Q$ such that $f+Q\not= 0$ on $b\D$. 
\vskip 2mm \rm\noindent Note that the only if part is a consequence of the argument principle.

One can view $f+Q$ above as $Pf+Q$ with $P\equiv 1$. We believe that an analogous 
theorem holds for $Q\equiv 1$:
\vskip 2mm
\noindent\bf CONJECTURE 1.2\ \ \it Let $f$ be a continuous function on $b\D$ such that
$$
W(Pf+1)\geq 0
\eqno (1.1)
$$
whenever $P$ is a polynomial such that $Pf+1\not= 0 $ on $b\D$. Then $f$ extends 
holomorphically through $\D $.
\vskip 2mm \rm 
\noindent The present note is the result of an unsuccessful attempt to prove 
this conjecture. In the paper we prove the conjecture for sufficiently smooth functions 
with finitely many zeros of finite order. In particular, the conjecture holds 
for real analytic functions. 
\vskip 4mm
\bf 2.\ Functions with no zeros
\vskip 2mm \rm
Suppose that the function $f$ has no zero. In this case (1.1) implies that $W(f)\geq 0$. Indeed,
there is an $\varepsilon>0$ such that $W(f)= W(f+\eta)$ for all $\eta,\ |\eta|<\varepsilon$. 
Choosing for $P$ a constant $c, \ |c|>1/\varepsilon $,\  (1.1) implies that $W(f)\geq 0$. 
Note that (1.1) implies that 
$$ W(f)+W\biggl(P+{1\over f}\biggr) = W (Pf+1)\geq 0
$$ so
$$
W(P+1/f)\geq -W(f) 
$$
for every polynomial $P$ such that $P+1/f\not= 0$ on $b\D$. 
If $W(f)=0$ then Theorem 1.1 implies that $1/f$ extends holomorphically through $\D$ and since $W(f)=0$ the argument principle shows that this 
holomorphic extension has no zero on $\D$ which gives
\vskip 2mm
\noindent\bf PROPOSITION 2.1\ \it Let $f$ be a continuous function 
on $b\D$ which has no zero and which satisfies $W(f)=0$. If $W(Pf+1)\geq 0 $ 
whenever $P$ is a polynomial such that 
$Pf+1\not= 0 $ on $b\D$ then $f$ extends holomorphically through $\D$. \rm
\vskip 2mm
\noindent Now, let $W(f)=N>0$. Then $W(Pf+1)=W(f) + W(P+1/f) \geq 0$ so
$
W(1/f + P)\geq -N
$
for every polynomial $P$ such that  $1/f+P\not= 0$ on $b\D$. A recent theorem of M.\ 
Raghupathi and M.\ Yattselev [RY, Th.2] applies to show that if $f$ is $\alpha$-H\" older 
continuous with $\alpha > 1/2$ then $1/f$ has a meromorphic extension through $\D$ which has 
at most $N$ poles, 
counting multiplicity. So in this case there are a function $H$ in the disc algebra and a 
polynomial $R$ of degree not exceeding $N$, 
with all zeros contained in $\D$, such that
$$
{1\over{f(z)}} = {{H(z)}\over {R(z)}} \ \ \ (z\in b\D ).
$$
Since $W(1/f)=-N$ and since $deg R \leq N$ it follows 
by the argument principle that $R$ has exactly $N$ zeros in $\D$, counting multiplicity, 
and $H$ has no zero on $\overline \D$. It follows 
that $f=R/H$ extends holomorphically through $\D$. This 
proves
\vskip 2mm
\noindent \bf PROPOSITION 2.2\ \ \it Let $f$ be an $\alpha$-H\" older continuous 
function on $b\D$ 
with $\alpha >1/2$ which has no zero on $b\D$. If \ $W(Pf+1)\geq 0$ whenever $P$ 
is a polynomial such that
$Pf+1\not=0$ on $b\D$ 
then $f$ extends holomorphically through $\Delta $. \rm
\vskip 4mm
\bf 3.\ Functions with finitely many zeros.
\rm
\vskip 2mm
The reasoning in Section 2 is no longer possible if $f$ has zeros on $b\D$. 
Suppose that $f$ has the form 
$$
f(z) = (z-b_1)^{m_1}(z-b_2)^{m_2}\cdots (z-b_n)^{m_n} g(z) \ \ \ (z\in b\D)
$$
where $b_i\in b\D ,\ 1\leq i\leq n$,\ $b_i\not=b_j\ (i\not=j)$   and where $g$ is a continuous function with no
zeros. (In particular, this holds when $f$ is real analytic). Assume that $W(g)=N$. Then
(1.1) implies
that 
$$
W(1/g + (z-b_1)^{m_1}\cdots (z-b_n)^{m_n}P)\geq -N .
$$
If $m_1=\cdots =m_n=0$ as in the preceding section then, if $N\geq 0$, [RY, Th.2] implies 
that $1/g$ has a meromorphic extension through $\Delta $ with at most $N$ poles. So 
the relevant question now is whether the same is true in general:
\vskip 2mm
\noindent \bf QUESTION 3.1 \it Let $b_
1, b_2, \cdots , b_n\in b\D,\ b_i\not=b_j$ if $i\not=j$, 
and let $m_1,\cdots, m_n\in\N$. Let $\cP$ be the family of all polynomials $Q$ of the form
$$
Q(z)=(z-b_1)^{m_1}\cdots (z-b_n)^{m_n}p(z)
$$
where $p$ is a polynomial, and let $J\in \N\cup\{ 0\}$. Suppose that $f$
is a continuous function on $b\D$ such that
$$
f(b_j)\not=0\ \ (1\leq j\leq n)
\eqno (3.1)
$$
and such that 
$
W(f+Q)\geq -J
$ for each $Q\in \cP$ such that $f+Q\not=0$ on $b\D$. Must $f$ extend meromorphically through $\D$ with the 
extension having at most $J$ poles, counting multiplicity? \rm
\vskip 2mm
\noindent Note that one has to assume (3.1) since otherwise $W(f+Q)$ is 
undefined for every $Q\in\cP$. If $m_1=\cdots = m_n=0$ and $J=0$ then
the positive answer 
is provided by Theorem 1.1. If $m_1=\cdots = m_n=0$ and $J\geq 1$ and 
if $f$ is $\alpha$-H\" older continuous 
with $\alpha > 1/2$ then the answer is positive by [RY,Th.2] which was 
 proved by using the theorem on rigid interpolaton:
\vskip 2mm
\noindent\bf THEOREM 3.2\ [RY, Th.5] \ \it Suppose that $g$ is a holomorphic function on $\D$ 
and let $N\in\N$. Suppose 
that for every nonnegative integer $n$ and for every polynomial $p$ of degree not exceeding 
$n$ we have 
$$
Z(z^ng+p)\leq N+n .
\eqno (3.2)
$$
Then $g$ is 
a quotient of polynomials of degree not exceeding $N$. \rm
\vskip 2mm
Note that if $g$ is a quotient of polynomials of degree not exceeding $N$ 
then (3.2) holds for every polynomial $p$ of degree not exceeding $n$. 
In the present paper we use Theorem 3.2 as an essential tool.
\vskip 4mm
\bf 4.\ On Question 3.1\rm
\vskip 2mm
\noindent\bf THEOREM 4.1\ \it Let $b_j\in b\D, \ 1\leq j\leq n, \ b_i\not=b_j\  (i\not=j)$,
and 
let $m_j\in\N, \ 1\leq j\leq n$.\ Let $\cP$ be the family of 
all polynomials $Q$ of the form 
$$ 
Q(z) = (z-b_1)^{m_1}(z-b_2)^{m_2}\cdots (z-b_n)^{m_n}p(z)
$$ 
where $p$ is a polynomial. Let $N= m_1+m_2+\cdots +m_n$ . Suppose that
$f\in \cC^{N+1} (b\D )$ is such that $f(b_j)\not= 0\ (1\leq j\leq n)$,
let $J$ be a nonnegative integer and assume that 
$$
W(f+Q)\geq -J\hbox{\ \ for every\ \ }Q\in\cP\hbox{\  such that\ \ } f+Q\not=0\hbox{\ on\ } b\D .
$$ 
Then $f$ has a meromorphic extension through $\D$ having at most 
$J$ poles in $\D$, counting multiplicity. \rm
\vskip 2mm
We shall rewrite the condition in Theorem 4.1 in a 
slightly different form. Let $a_j\in b\D,\ 1\leq j\leq N$, and
let $\cS$ be the family of all polynomials $Q$ of the form
$$
Q(z) = (z-a_1)\cdots (z-a_N) p(z)
$$
where $p$ is a polynomial. Note that we do not require 
that $a_i\not= a_j$ if $i\not= j$. 
To prove Theorem 4.1 it will be enough to prove the following 
$$\left.
\eqalign{&\hbox{Suppose that\ } f \hbox{\ is a function of class\ } \cC^{N+1} \hbox {\ on\ }
b\D \hbox{\ such that\ } 
f(a_j)\not= 0,  \cr &1\leq j\leq N, \hbox{\ and such that \ }W(f+Q)\geq -J \hbox{\ whenever\ } 
Q\in\cS\hbox{\ and\ } f+Q\not=0\cr &\hbox{on\ }b\D. \hbox{\ Then\ } f \hbox{\ has 
a meromorphic extension through 
\ } \D \hbox{\ having at most\ }J\cr &\hbox{poles in\ } \D \hbox{\ counting multiplicity. \ 
}\cr}\right\} \eqno (4.1)
$$
\vskip 1mm
\noindent \bf PROPOSITION 4.2\ \it A continuous function $f$ on $b\D$
satisfies 
$$
W(f+(z-a_1)\cdots (z-a_N)g)\geq -J
\eqno (4.2)
$$
for every polynomial $g$ such that
$$ f + (z-a_1)\cdots (z-a_N)g\not= 0 \hbox{\ on\ } b\D. 
\eqno (4.3)
$$
if and only if $f$ satisfies \rm(4.2)\it \ for every function $g\in A(\D )$ which satisfies 
\rm(4.3).
\vskip 1mm
\noindent\bf Proof. \rm If for some $g\in A(\D)$ satisfying (4.3) we have  
$W( f + (z-a_1)\cdots (z-a_n)g)\leq -J-1$ then this holds for all sufficiently small perturbations of $g$.
In particular, it holds for some polynomial $g$, contradicting (4.2), completing the proof. 
\vskip 2mm
\noindent\bf LEMMA 4.3\ \it Let $I\subset b\D$ be an arc centered at $a$ and let $f\in\cC ^{n+1}
(I)$. There are a polynomial $p $ of degree not exceeding $n-1$ and a
function $h\in \cC^1 (I)$ such that
$f(z) = p(z) + (z-a)^n h(z)\ \ (z\in I)$. 
\vskip 1mm \noindent \bf Proof. \rm 
Write $a=e^{it_0}$ and let $J$ be a segment on $\R$ centered at $t_0$. Then
$$
f(e^{it})=f(e^{it_0}) + c_1(t-t_0)+\cdots + c_{n-1}(t-t_0)^{n-1} + (t-t_0)^{n}g(t)\ \ (t\in J)
\eqno (4.4)
$$
where $g\in\cC^1 (J)$. For all $t$ close to $t_0,\ t\not= t_0$ we have
$$
{{t-t_0}\over{e^{it}-e^{it_0}}}= a_0+a_1 (e^{it}-e^{it_0})+a_2(e^{it}-e^{it_0})^2+\cdots
\eqno (4.5)
$$
where the series converges for $t$ near $t_0$. Now, by (4.4) 
$$\eqalign{
f(e^{it}) = f(e^{it_0}) + c_1 \biggl( {{t-t_0}\over{e^{it}-e^{it_0}}}\biggr)(e^{it}-e^{it_0})
+\cdots +c_{n-1}
&\biggl({{t-t_0}\over{e^{it}-e^{it_0}}}\biggr)^{n-1}({e^{it}-e^{it_0}})^{n-1}+\cr
&
+\biggl({{t-t_0}\over{e^{it}-e^{it_0}}}\biggr)^{n}({e^{it}-e^{it_0}})^{n}g(t)\cr}
$$
which, by (4.5) implies that
$$
\eqalign{f(e^{it})= f(e^{it_0})+&c_1\biggl[\sum_{j=0}^\infty
a_j ({e^{it}-e^{it_0}})^j\biggr]({e^{it}-e^{it_0}})+\cr
+&c_2\biggl[\sum_{j=0}^\infty
a_j ({e^{it}-e^{it_0}})^j\biggr]^2({e^{it}-e^{it_0}})^2 + \cr
&+ \cdots +\cr
+&c_{n-1}\biggl[\sum_{j=0}^\infty
a_j ({e^{it}-e^{it_0}})^j\biggr]^{n-1}({e^{it}-e^{it_0}})^{n-1}+\cr
+&\biggl[\sum_{j=0}^\infty
a_j ({e^{it}-e^{it_0}})^j\biggr]^n({e^{it}-e^{it_0}})^n g(t)\cr }
$$
Computing the powers and rearranging we get
$$f(e^{it}) = f(e^{it_0}) + b_1({e^{it}-e^{it_0}})+\cdots + b_{n-1}({e^{it}-e^{it_0}})^{n-1} +
({e^{it}-e^{it_0}})^n \bigl[ g(t)+w(t)\bigr]
$$
where, as a sum of a convergent power series $w$ is real analytic so that $g+w=h$ is of class 
$\cC^1$.
The proof is complete.
\vskip 2mm 
\noindent \bf LEMMA 4.4\ \it Let $f$ be a function of class
$\cC ^{N+1} $ on $b\D$ and let $b_j\in b\D,\ 1\leq j\leq n,\ b_i\not= b_j\ (j\not= i)$. 
Let $m_1,\cdots m_n$ be positive integers such that $m_1+\cdots m_n=N$ and let
$$
\eqalign{
&a_k=b_1\ \ (1\leq k\leq m_1) \cr
&a_k=b_2\ \ (m_1+1\leq k\leq m_1+m_2) \cr
&\cdots \cr
& a_k=b_n\ \ (m_1+\cdots m_{n-1}+1\leq k\leq m_1+\cdots m_n .)\cr}
$$
There are constants $A_k,\ 0\leq k\leq N-1$, and a function of class $\cC^1 $ on $b\D$ such that
$$
\eqalign{ &f(z) = A_0+A_1 (z-a_1)+ A_2(z-a_1)
(z-a_2)+\cr &\ \ \ +\cdots + A_{N-1}(z-a_1)\cdots (z-a_{N-1})
+ g(z)
(z-a_1)\cdots (z-a_N).\cr}
$$
\vskip 1mm
\noindent \bf Proof.\ \rm By Lemma 4.3 we have
$$
f(z)= B_{10}+ B_{1 1}(z-b_1)+\cdots +B_{1, m_1-1}(z-b_1)^{m_1-1} + g_1(z)(z-b_1)^{m_1}\ \ (z\in b\D)
$$ 
where $g_1$ is of class $\cC^1$ on $b\D$ and of class $\cC^{N+1}$ on $b\D\setminus \{ b_1\} $. We repeat the procedure to write
$$
g_1(z)= B_{20}+ B_{2 1}(z-b_2)+\cdots +B_{2, m_2-1}(z-b_2)^{m_2-1} + g_2(z)(z-b_2)^{m_2}\ \ (z\in b\D)
$$
where the function $g_2$ is of class $\cC^1$ on $b\D$ and of class $\cC^{N+1}$ on
$b\D \setminus \{ b_1, b_2\} $. Repeating this procedure we get the functions $g_2, g_3,
\cdots , g_n$, all of class
$\cC^1 $ on $b\D$, 
such that
$$
g_{n-1}(z)= B_{n0}+ B_{n 1}(z-b_n)+\cdots +B_{n, m_n-1}(z-b_n)^{m_n-1} + g_n(z)(z-b_n)^{m_n} .
$$
Putting $g=g_n$ and substituting the expression for $g_{n-1}$ into the expression for $g_{n-2}$ and so on 
we get the result with
$$
\eqalign{ 
&A_0=B_{10},\cdots , A_{m_1-1}=B_{1,m_1-1}, \cr
&A_{m_1}=B_{20},\ A_{m_1 + 1}=B_{2 1} ,\cdots , A_{m_2-1}= B_{2, m_2-1},\cr
&\cdots \cr}
$$
which completes the proof.
\vskip 4mm
\bf 5.\ Proof of Theorem 4.1
\vskip 2mm
\rm As already mentioned we have to prove (4.1). So suppose that 
$f\in \cC^{N+1}(b\D)$ satisfies $f(a_j)\not=0\ 
(1 \leq j\leq N)$ and satisfies (4.2) for every polynomial $g$ satisfying (4.3). By Proposition 4.2
$f$ satisfies (4.2) for 
every $g$ in the disc algebra that satisfies (4.3). By Lemma 4.4 there
are numbers $A_0, A_1,\cdots  A_{N-1}$  and 
a function $g$ of class $\cC ^1$ on $b\D$ such that if
$$
D(z)=A_0+A_1(z-a_1)+\cdots + A_{N-1}(z-a_1)\cdots (z-a_{N-1})
$$
then
$$
f(z) = D(z) +
(z-a_1)\cdots (z-a_N) g(z)\ \ (z\in b\D) 
$$
which implies that 
$$
W\bigl(D
+(z-a_1)\cdots (z-a_N)g + (z-a_1)\cdots (z-a_N)P\bigr) \geq -J
$$ 
for every $P$ in the disc algebra such that the the expression in the parenthesis is different 
from $0$ on $b\D$.  
Since $g$ is of class $\cC^1$ we can write
$$
g(z)= F(z)+\overline {G(z)}\ \ \ (z\in b\D )
$$
where $F$ and $G$ are in the disc algebra
with boundary values of class $H^\alpha$ for every $\alpha <1$ 
which implies that 
$$
W\bigl(D + (z-a_1)\cdots (z-a_N)(F+\overline G+P)\bigr)\geq -J
$$
for every $P$ in the disc algebra such that the expression in the parenthesis is different 
from $0$ on $b\D$ so that 
$$
W\bigl(D +
 (z-a_1)\cdots (z-a_N)(\overline G+P)\bigr)\geq -J
\eqno (5.1)
$$
whenever $P$ in the disc algebra is such that
$$
D 
+(z-a_1)\cdots (z-a_N)(\overline G+P))\not= 0 \hbox{\ \ on\ \  }b\D .
\eqno (5.2)
$$
By our assumption we have $f(a_j) \not= 0\ (1\leq j\leq N)$ which implies that
$$
D(a_j)\not= 0\ \ \ (1\leq j\leq N).
\eqno (5.3)
$$
Recall that (5.1) holds whenever $P$ is in the disc algebra and is such that (5.2) holds. 
In particular, (5.1) holds whenever $P$ is a polynomial satisfying (5.2). 
Conjugating (5.1) we get
$$
W\bigl(\overline{D}
 +(\overline z-\overline{a_1})\cdots (\overline z-\overline{a_N})( G+\overline P)\bigr)\leq J
$$
which, multiplying the expression in the parenthesis with 
$z^N$ gives
$$
W\bigl(z^N\overline{D} + 
(1-\overline{a_1}z)\cdots (1-\overline{a_{N}}z)(G+\overline P)\bigr)\leq N+J 
$$
for every polynomial $P$ such that the expression in the parenthesis does not vanish on $b\D$.
In particular, if $A$ is the polynomial such that
 
$$ 
A(z)=z^N\overline{D(z)}\ \ (z\in b\D ) \hbox{\ \ and if\ } \ \ 
B(z)= (1-\overline{a_1}z)\cdots (1-\overline{a_{N}}z)
$$
then, the degree of $A$ does not exceed $N$ and we have 
 $$ 
 W(A+BG+B\overline P)  \leq N+J
 $$ 
  for every polynomial  $P$ such that $A+BG+B\overline P\not=0$ on $b\D$. On $b\D$ we have
  $\overline {P(z)} = p(z)/z^m$ where $m\in\N\cup\{ 0\}$ and where $p$ is a polynomial of
  degree not exceeding $m$ so it follows that 
  $$
  W(z^m(A+BG)+Bp_m)\leq N+J+m
  $$
  whenever $m\in\N\cup\{ 0\}$ and $p_m$ is a polynomial of degree not exceeding $m$ such that
  $$
  z^m(A+BG)+Bp_m\not = 0 \hbox{\ on\ } b\D
  \eqno (5.4)
  $$
The argument principle implies that 
$$
Z\biggl(z^m\biggl({A\over B}+G\biggr)+p_m\biggr)\leq N+J+m
\eqno (5.5)
$$
for every $m\in\N\cup\{ 0\}$ and for every polynomial $p_m$ of degree not exceeding $m$ such that (5.4) holds. 
Now, by (5.3) we have  $A(a_j)\not= 0\ (1\leq j\leq N)$. Since $B(a_j)=0 \ (1\leq j\leq N)$, \ 
(5.4) holds if and only if 
$$
z^m\biggl( {{A(z)}\over {B(z)}}+G(z)\biggr)+p_m(z) \not=0 \ \ (z\in b\D,\ z\not=a_j ,\ 
1\leq j\leq N)
\eqno (5.6)
$$
It follows that (5.5) holds for every polynomial $p_m$ of degree not exceeding $m$,
without condition (5.6). 
Indeed if 
$$
Z\biggl( z^m\biggl({A\over B}+G\biggr)+q_m\biggr)\geq N+J+m+1
\eqno (5.7)
$$
for some $m\in\N\cup\{ 0\}$ and for some polynomial $q_m$ of degree not exceeding $m$ then, by the argument principle, the same holds for $q_m$ replaced 
by $q_m +\eta$ for all sufficiently small $\eta$. However, since $G$ is $\alpha$-H\" older smooth with $\alpha >1/2$ 
the same holds for the 
function $z\mapsto z^m\bigl(A(z)/B(z)+G(z)\bigr) + q_m(z)$ which implies that 
the set
$$
\biggl\{ z^m\biggl( {{A(z)}\over {B(z)}}+G(z)\biggr)+q_m(z)
\colon z\in b\D,\ z\not= a_j,\ 1\leq j\leq N \biggr\}
$$ has planar measure zero so there are arbitrarily small $\eta$ such that (5.7) holds 
for $q_m$ replaced by $p_m=q_m+\eta$
where $p_m$ satisfies (5.6). This proves that (5.5) holds for  any $m\in\N\cup\{ 0\}$ 
and for any polynomial $p_m$ of degree not exceeding $m$. Theorem 3.2 now 
applies to show that there are polynomials $P, Q$ of degree not 
exceeding $N+J$ without common factors such that
$$
{{A(z)}\over{ B(z)}} + G(z) = {{P(z)}\over {Q(z)}}\ \ (z\in b\D ) .
$$
Recall that 
$$ 
B(z) = {1\over{a_1a_2\cdots a_n}}(a_1-z)(a_2-z)\cdots (a_N -z)
$$
and that $A(a_j)\not= 0\ (1\leq j\leq N)$ . We have
$G= P/Q-A/B$. $G$ is continuous on $\overline \D$ so if 
a factor $(\alpha-z)^k $ occurs in $B$ then $Q$ has to 
be divisible by $(\alpha-z)^k $. In fact, if $Q$ contained only 
$\ell$ factors $(\alpha -z)$ with $\ell<k$ then we would get
$$
(z-\alpha )^\ell G = {P\over{Q_1}} - {A\over{(z-\alpha )^{k-\ell}}B_1}
$$
where $Q_1$ is a polynomial, $Q_1(\alpha )\not=0$, and where $B_1$ is a polynomial. Since the left side  is continuous at $\alpha$, since $Q_1(\alpha)\not=0$ and since 
$A(\alpha)\not=0$ this is not possible. Thus $Q=RB$ where R is a polynomial of degree not exceeding $J$. 
It follows that 
$$
G={P\over Q} - {A\over B} = {{P-RA}\over {RB}}.
$$
Since $G$ is continuous on $\overline\D$ it follows that $P-RA$ must be 
divisible by $B$. Since $deg A\leq N,\ deg R\leq J, \ deg P\leq N+J$ it 
follows that $P-RA$ is a polynomial of degree not exceeding $N+J$. Thus, on $b\D$, \ $G$ is 
a quotient of two polynomials 
of degree not exceeding $J$ which implies that the same holds for $\overline G$. Thus, 
$g= F+\overline G$ extends meromorphically through $\D$ with 
the number of poles not exceeding $J$ and so the same holds for 
$$
f = A_0+A_1(z-a_1)+\cdots + A_{N-1}(z-a_1)\cdots (z-a_{N-1}) + (z-a_1)\cdots (z-a_N)g .
$$
The proof is complete.
\vskip 2mm
\noindent \bf REMARK \ \rm Another look at the proof of Theorem 4.1 
shows that it is enough to assume that
$f\in \cC^{L+1}(b\D)$ where $L= \max \{ m_1,m_2,\cdots, m_n\}$.
\vskip 4mm
\bf 6.\ On Conjecture 1.2\rm
\vskip 2mm
We prove somewhat more general result than Conjecture 1.2 for 
sufficiently smooth functions $f$ with finitely many zeros on $b\D$ of finite order:
\vskip 2mm
\noindent\bf THEOREM 6.1 \ \it Suppose that $f$ is of class $\cC^\infty$ on $b\D$ with at most 
finitely many zeros of finite order and let $J\in \N\cup\{ 0\}$. Then
$$
W(Pf+1)\geq -J
$$
for each polynomial $P$ such that $Pf+1\not= 0$ on $b\D$ if and only if $f$ 
extends meromorphically 
through $\D$ with the extension having at most $J$ poles, counting multiplicity. 
In particular, $f$ 
extends holomorphically through $\D$ if and only if $W(Pf+1)\geq 0$ for every 
polynomial $P$ such that $Pf+1\not=0$ on 
$b\D$. 
\vskip 2mm
\noindent\bf COROLLARY 6.2\ \it A real analytic function $f$ on $b\D$ extends 
meromorphically through $\D$ with at most $J$
poles in $\D $, counting multiplicity, if and only if $W(Pf+1)\geq -J$ for 
each polynomial $P$ such that $Pf+1\not= 0$ on 
$b\D$. In particular, $f$ extends holomorphically through $\D$ if and only 
if $W(Pf+1)\geq 0$  for each polynomial $P$ 
such that $Pf+1\not=0$ on $b\D$. \rm
\vskip 2mm
\noindent\bf Proof of Theorem 6.1. \rm Let $J$ be a nonnegative 
integer and let $f$ be a smooth function on $b\D$ that satisfies
$$
W(Pf+1)\geq -J
\eqno (6.1)
$$
for all polynomials $P$ such that $Pf+1\not= 0$ on $b\D$. This 
happens if and only if 
(6.1) holds for all functions $P$ in the disc algebra such that
$Pf+1\not=0$ on  $b\D$. Indeed,
if we have $W(Pf+1)\leq -J-1$ for some $P_0$ in the disc algebra 
then the same holds for all $P$
in the disc algebra sufficiently close to $P_0$. In particular, 
it holds for some polynomial $P$. We assume that $f$ has at most finitely
many zeros (of finite order) on
$b\D$ so $f=\Pi g$ where $\Pi (z) = (z-a_1)(z-a_2)\cdots (z-a_n)$ 
and $g$ is a smooth function 
on $b\D$ without zeros.  Now
(6.1) becomes
$$
W(P\Pi g+1)\geq -J
$$
which gives
$$
W\biggl( P\Pi + {1\over g} \biggr) \geq -J-W(g).
$$
Suppose first that $W(g)=N\geq -J$ so that $N+J\geq 0$ and so
$$
W\biggl( {1\over g} + \Pi P\biggr) \geq -(N+J)
$$
for every polynomial P such that $1/g + \Pi P\not=0$ on $b\D$.  By Theorem 4.1
$$
{1\over {g(z)}} = {H(z)\over {Q(z)}}\ \ (z\in b\D)
$$
where $H$ is in the disc algebra and $Q$ is a polynomial with at most $N+J$ zeros on $\D$. 
The argument principle now shows that $N=W(g)=W(Q)-
W(H) = Z(Q)-Z(H) \leq N+J -Z(H)$ which implies that $Z(H)\leq J$ which 
shows that $g$, and consequently $f=\Pi g$ has 
a meromorphic extension through $\D$ with at most $J$ poles, counting multiplicity. 

We now complete the proof by showing that $N+J<0$ is impossible. Assume that $-(N+J)\geq 1$. 
Since $g$ is smooth and $W(g)=N <0$ one can write 
$$
g(z) = F(z)\overline {G(z)}z^N\ \ (z\in b\D)
$$
where $F$ and $G$ are in the disc algebra with no zeros on $\overline \D$ and 
with smooth boundary values. We get
$$
W(P\Pi F \overline Gz^N +1)\geq -J
$$
whenever P in the disc algebra is such that the expression in the parenthesis is 
different from zero on $b\D$. Since $F$ has no zero on $\overline \D$ this happens if and only if
$$
W(P\Pi\overline G + z^{-N})\geq -J-N\geq 1
$$
and, since $G$ has no zero on $\overline \D$ we have
$$
W\biggl( {{z^{-N}}\over {\overline G}}+P\Pi\biggr)\geq 1
$$
whenever $P$ in the disc algebra is such that 
$$
{{z^{-N}}\over{\overline G}} + P\Pi \not= 0\hbox{\ on\ } b\D .
$$ 
By Theorem 4.1 it follows that $z^{-N}/\overline G$ extends holomorphically through $\D$ which is possible if 
and only if 
$$
{{z^{-N}}\over{\overline {G(z)}}} = Q(z)\ \ (z\in b\D)
$$
where $Q$ is a polynomial of degree not exceeding $-N$.
In particular,
$$
W(Q +P\Pi)\geq 1
\eqno (6.2)
$$
for all functions $P$ in the disc algebra such that $Q+P\Pi\not=0$ on $b\D $. However, 
one can choose $P$ in the disc algebra such that $Q+P\Pi = e^\Psi$ with $\Psi$ entire. To
do this one has to choose $\Psi $ in such a way that $( e^\Psi -Q)/\Pi$ is holomorphic. This 
is easy to do, see [G4]. 
The argument principle now shows that with this $P$, we have $W(F+P\Pi )=0$ so (6.2) fails. 
This shows that $N+J<0$ is impossible and completes the proof of Theorem 6.1.
\vskip 4mm
\bf 7. Generalizations of Theorem 3.1 \rm
\vskip 2mm
Suppose that $f$ is a continuous function and $R$ is a fixed polynomial. 
Suppose that
$f$ is a continuous function on $b\D$ that does not vanish at any zero of 
$R$ contained in $b\D$, 
such that
$$
W(f+Rp)\geq 0
$$
for every polynomial $p$ such that $f+Rp\not= 0 $ on $b\D $. 
Must $f$ extend holomorphically through $\D$? We know from Section 4 that 
the answer is positive provided that all zeros of $R$ are on $b\D$ and 
provided that $f$ is sufficiently smooth. 

For general R the answer is negative. To see this, 
let $f(z)=z/(z-1/2)\ (z\in b\D)$. If $p$ is a polynomial such 
that $f+zp\not= 0$ on $b\D$ then 
the argument principle implies that 
$$\eqalign{
W(f+zp) = W \biggl({z\over{z-1/2}} + zp\biggr) = 
W\biggl( {z\over{z-1/2}}\biggl( 1+(z-1/2)p\biggr)\biggr) =\cr
=W\biggl( {z\over{z-1/2}}\biggr) + W\bigl( 1+(z-1/2)p) \geq 0 \cr}
$$
yet $f$ does not extend holomorphically through $\D$. We now show that for sufficiently 
smooth functions the answer to the question above is positive provided that $R$ has no zero in 
$\D$. 

Let $\Pi_1$ be a product of $N\geq 0$ factors of the form $z-a, a\in\D$, let $\Pi_2 $ be a 
finite product of factors of the form $z-a,\ a\in b\D$,  and let $\Pi_3$ be a 
finite product of factors of the
form $z-a,\ a\in \C\setminus\overline\D$. Let $\Pi =\Pi_1\Pi_2\Pi_3$, let $J$ be a 
nonnegative integer 
and suppose that $f$ is a smooth function on 
$b\D$ such that $f$ does not vanish at any zero of $\Pi_2$ and such that
$$
W(f+\Pi p)\geq -J
\eqno (7.1)
$$ whenever $p$ is a polynomial such that $f+\Pi p\not= 0$ on $b\D$. 
We know that this happens if and only if (7.1) holds for each $p$ in the disc algebra such that
$f+\Pi p\not= 0$ on $b\D$. Now, since the zeros of $\Pi_3$ are in $\C\setminus\overline\D$ it 
follows that $p$ is in the disc algebra if and only $\Pi_3 p$ is in the disc algebra. 
It follows that 
(7.1) holds for every $p$ in the disc algebra 
such that $f+\Pi p\not= 0$ on $b\D$ if and only if 
$$
W(f+\Pi_1\Pi_2p)\geq -J
\eqno (7.2)
$$
for each $p$ in the disc algebra such that $f+\Pi_1\Pi_2p\not=0 $ on $b\D $. Now, (7.2)
implies that
$$
W\biggl( {f\over{\Pi_1}} + \Pi_2p\biggr)\geq -J-N
$$
whenever $p$ is a polynomial such that 
$f/\Pi_1+\Pi_2p\not=0$ on $b\D$. If $f$ is sufficiently smooth then 
Theorem 4.1 implies that $f/\Pi_1$ has a 
meromorphic extension through $\D $ which has at most $J+N$ poles in $\D $, 
counting multiplicity.
This proves
\vskip 2mm
\noindent\bf THEOREM 7.1\ \it Let $\Pi=\Pi_1\Pi_2\Pi_3$ where $\Pi_1$ is a 
product of $N$ factors 
of the form $z-a,\ a\in\D$, where $\Pi_2$ is a finite product of factors 
of the form $z-a,\ a\in b\D$, and where $\Pi_3$ is a finite product of factors 
of the form $z-a,\ a\in\C\setminus \overline\D$. Assume that $f\in \cC^\infty (b\D)$ 
vanishes at no zero of 
$\Pi_2$ and assume that $J$ is a nonnegative integer. Then $f$ satisfies
$$ 
W(f+\Pi p)\geq -J
$$
for every polynomial $p$ such that $f+\Pi p\not=0$ on $b\D $ if and only if $f/\Pi_1 $ 
has a meromorphic 
extension through $\D$ which has at most $N+J$ poles in $\D$, counting multiplicity.
\vskip 4mm
\bf 8. Remarks \rm
\vskip 2mm
If $D$ is a bounded domain in $\C$ we denote
by $A(D)$ the algebra of all continuous functions on
$\overline D$ which extend holomorphically through $D$. 
Theorem 1.1 has been generalized to
\vskip 2mm
\noindent\bf THEOREM 8.1\rm \ [G5]\ \it Let $D\subset \C$ be a bounded domain whose
boundary consists of 
a finite number of pairwise disjoint simple closed curves. Let $J$ be a nonnegative 
integer. Then 
$W(Pf+Q)\geq -J$ for each $P, Q$ in $A(D)$ such that $Pf+Q\not= 0$ on $bD$ if and 
only if $f$ has a meromorphic extension through $D$ with at most $J$ poles 
counting multiplicity. \rm
\vskip 2mm
\noindent If $J=0$, that is, if we speak of holomorphic extendibility, then 
one can take $P\equiv 1$ [G2]. It remains an open question whether one can take 
$P\equiv 1$ in general. 
Ragupathi and Yattselev [RY] made progress by proving that one can
take $P\equiv 1$ in the case when $D=\D $ and $f$ is $\alpha$-H\" older continuous
with $\alpha>1/2. $ 
Conjecture 1.2 deals with the open question whether one can take $Q\equiv 1$ in 
Theorem 8.1. 
We conclude by mentioning a related result which holds for all continuous 
functions:
\vskip 2mm
\noindent \bf THEOREM 8.2\ \it Let $f$ be a continuous function on $b\D $ and assume that
$$
W\bigl( P(f+c)+1)\geq 0
$$
whenever $c$ is a constant and $P$ is a polynomial such that $P(f+c)+1\not= 0$ on $b\D$.
Then $f$ extends holomorphically through $\D$. \rm
\vskip 1mm \noindent\bf Proof.\ \rm Observe that by choosing $c$ large enough,
$W(f+c)=0$ and so, 
by Proposition 2.1, $f+c$ extends holomorphically through $\D$ and so does $f$.
\vskip 5mm
\rm This paper was supported in part by the ministry of Higher Education, 
Science and Technology of Slovenia through the research
  program Analysis and Geometry, Contract No. P1-02091 (B).
\vfill
\eject
\centerline{\bf REFERENCES}\rm 
\vskip 3mm

\noindent [G1]\ J.\ Globevnik:\ Holomorphic extendibility and the argument principle.

\noindent Contemp.\ Math.\  Vol.\ 382 (2005) 171-175
\vskip 2mm
\noindent [G2]\  J.\ Globevnik:\ The argument principle and holomorphic extendibility.

\noindent Journ.\ d'Analyse.\ Math.\ 94 (2004) 385-395
\vskip 2mm
\noindent [G3]\ J.\ Globevnik:\ The argument principle and holomorphic extendibility 
to finite Riemann surfaces.

\noindent  Math.\ Z.\ 253 (2006) 219-225
\vskip 2mm
\noindent [G4]\ J.\ Globevnik:\ Meromorphic extendibility and the argument principle.

\noindent Publ.\ Mat.\ 52 (2008) 171-188
\vskip 2mm
\noindent [G5]\ J.\ Globevnik;\ On meromorphic extendibility.

\noindent J.\ Math.\ Anal.\ Appl.\ 351 (2009) 285-290
\vskip 2mm
\noindent [K]\ D.\ Khavinson:\ A note on a theorem of J. Globevnik

\noindent Contemp.\ Math.\  Vol.\ 382 (2005) 227-228
\vskip 2mm
\noindent [RY]\ M.\ Raghupathi and M.L.Yattselev:
\ Meromorphic extendibility and rigidity of interpolation.

\noindent J.\ Math.\ Anal.\ Appl.\ 377 (2011) 828-833
\vskip 2mm
\noindent [S]\ E.\ L.\ Stout: \it The theory of uniform algebras.\rm 

\noindent Bogden and Quigley, Tarrytown-on-Hudson, NY, 1971
\vskip 2mm
\noindent [T]\ M.\ Tsuji:\ \it Potential theory in modern function theory.\rm

\noindent Maruzen, Tokyo, 1959
\vskip 2mm
\noindent [Z] A.\ Zygmund: \it Trigonometric series. \rm

\noindent Cambridge University Press, Cambridge, New York, 1959

\vskip 8mm
\noindent Institute of Mathematics, Physics and Mechanics

\noindent University of Ljubljana, Ljubljana, Slovenia

\noindent josip.globevnik@fmf.uni-lj.si

\bye